\documentclass[12pt]{article}
\usepackage{latexsym}
\usepackage{amssymb}
 
\newcommand{\proof}{{\noindent \bf Proof. }}
\newtheorem{thm}{Theorem}

\newtheorem{cor}{Corollary}
\newtheorem{prop}{Proposition}

\newcommand{\x}{{\bf x}}
\newcommand{\y}{{\bf y}}

\newcommand{\E}{{\mathbb E}}
\newcommand{\N}{{\mathbb N}}

\newcommand{\B}{{\cal B}}

\newcommand{\M}{{\cal M}}

\date{}

\begin{document}
\begin{titlepage}
\title{\bf Skewincidence}
{\author{{\bf  G\'erard Cohen}
\\{\tt cohen@enst.fr}
\\ENST
\\ FRANCE
\and{\bf Emanuela Fachini}
\\{\tt fachini@di.uniroma1.it}
\\''La Sapienza'' University of Rome
\\ ITALY
\and{\bf J\'anos K\"orner}
\thanks{Department of Computer Science, University of Rome, La Sapienza, 
via Salaria 113, 00198 Rome, ITALY}
\\{\tt korner@di.uniroma1.it}
\\''La Sapienza'' University of Rome
\\ ITALY}} 

\maketitle
\begin{abstract}
We introduce a new class of problems lying halfway between questions about graph capacity and intersection.  We say that two binary sequences $\x$ and 
$\y$ of the same length have a skewincidence if there is a coordinate $i$ for which $x_i=y_{i+1}=1$ or vice versa. We give rather sharp bounds on the maximum number of binary sequences of length $n$ any pair of which has a skewincidence.
\end{abstract}
\end{titlepage}
 
\section{Introduction}

We consider binary relations of strings of some fixed finite length $n$ from a finite alphabet (or strings representing the linear orders of $[n]$). We are interested in the maximum number of strings any two of which are in the given relation. Most problems of this kind belong to one of two well--investigated classes of opposite nature. 

{\it Intersection problems} have been studied in extremal combinatorics. The first of these goes back to the seminal paper of Erd\H os, Ko and Rado \cite{EKR}. These authors say that the binary strings $\x=x_1x_2\dots x_n$ and $\y=y_1y_2\dots y_n$ intersect if for some coordinate $i$ they have $x_i=y_i=1$. 
They then determine the maximum number of pairwise intersecting binary strings of length $n$ and weight $k$; here the weight of a string is its number of 1's. (In other words, they determine the largest stable sets in the Kneser graph whose vertices are the elements of ${[n] \choose k}$.) They show that any optimal configuration has the same structure; it consists of all those strings that have a 1 in the same fixed position. In other words, these seuquences have a 
fixed projection on some coordinate. Such a structure is often called a {\it kernel structure} and it is the natural candidate solution for all the intersection problems. The reason for this seems to be the fact that the relation underlying the problem is a {\it similarity relation}. We will say that a binary relation for strings of the same length is a similarity relation if it is reflexive and locally verifiable, meaning that if some projections of two strings are in this relation 
then this implies that so are the strings themselves. For more on this, we refer to \cite{CFGS} and \cite{EFP}. 

{\it Capacity problems} originate in the fundamental paper of Claude Shannon \cite{Sh} and come from information theory. We will say that a binary relation for strings of the same length is a difference relation if the relation is irreflexive and locally verifiable \cite{FK}. For easy reference, we will say that two sequences are very different if they are in the given difference relation. For a fixed length, one is interested, as before, in the maximum number of pairwise very different sequences. The classical example comes from Shannon and has been generalized in a series of papers; for more on this we refer to 
\cite{CKS} and the survey \cite{KO}. Unlike for intersection problems, here there is no natural conjecture for the optimal constructions and most problems of this kind 
remain wide open.

Both of these groups of problems have been generalized in recent work to permutations of $[n]$. For intersection problems on permutations we refer to 
\cite{EFP}. Capacity problems for permutations have been introduced in \cite{KM}; for further developments we cite \cite{BCFFKST}.
In order to introduce our new problems it will be interesting to recall the first capacity problem on permutations from \cite{KM}. 
We call two permutations of $[n]$ colliding if they map some $i\in [n]$ into two consecutive integers. Let us denote by 
$T(n)$ the maximum cardinality of a set of pairwise colliding permutations of $[n]$. K\"orner and Malvenuto \cite{KM} conjecture that 
$T(n)$ equals the middle binomial coefficent ${n \choose {\lfloor \frac{n}{2} \rfloor}}$. This conjecture is still open; for the best bounds 
we refer to \cite{BCFFKST}.

Our starting point in the present work is the problem about colliding permutations. We note that if two permutations, $\rho$ and $\sigma$ 
are colliding, then their inverses are skewincident. In fact, the collision relation means that for some $j \in [n]$ we have 
$$|\rho(j)-\sigma(j)|=1$$ Suppose without loss of generality that $\sigma(j)=\rho(j)+1$. Denoting $i:=\rho(j)$ we have $\rho^{-1}(i)=\sigma^{-1}(i+1)=j$
meaning that there is a skewincidence between the strings describing the two permutations; we find in them the same symbol $j$ in adjacent positions. 
The resulting relation of coincidence is irreflexive for permutations. For sequences with repetitions such as long strings from a finite alphabet the 
analogous relation is not irreflexive any more. In fact, it is neither reflexive nor irreflexive and as our initial findings show the optimal solution 
has a somewhat unusual behaviour. Our results are asymptotic. Logarithms and exponentials are to the base 2.

\section{Results}

Let us fix a natural number $n$ and consider the set $\{0,1\}^n$ of the binary strings of length $n$. We say that the sequences $\x \in \{0,1\}^n$ and 
$\y \in \{0,1\}^n$ have a {\it skew coincidence} (abbreviated as skewincidence) if for some coordinate $i \in [n-1]$ we have either $x_i=y_{i+1}=1$ or $x_{i+1}=y_i=1.$ Let us 
denote by $M(n)$ the maximum number of binary strings of length $n$ any two of which have a skew coincidence. We have the following result

\begin{thm}\label{thm:main}
$$2^n-2^{0.96n}\leq M(n)\leq 2^n-2^{0.69n}$$
for $n$ sufficiently large.
\end{thm}

This implies

\begin{cor}
\label{cor:mainc}
$$\lim_{n\rightarrow \infty}\frac{M(n)}{2^n}=1$$
\end{cor}

\proof

To prove the upper bound, let us consider the set $\sf{F}_n\subseteq \{0,1\}^n$ of those binary sequences that do not contain a 1 in consecutive positions. 
It is well--known that
$$\vert {\sf F}_n \vert=f_n$$
where $f_1=2$, $f_2=3$, $f_n=f_{n-1}+f_{n-2}$  meaning that $\{f_n\}_{n=1}^{\infty}$ is the standard Fibonacci sequence. Given two binary sequences $\x$ and $\y$ we say that $\x \leq \y$ 
if $\x=x_1x_2 \dots x_n$, $\y=y_1y_2 \dots y_n$, and $x_i \leq y_i$ for every $i \in [n]$. We say that $\x$ and $\y$ are comparable if $\x \leq \y$ or vice versa.
Consider now a set $\sf{B}$ of pairwise skewincident binary strings 
from $\{0,1\}^n$. It is obvious that if two strings belong to the intersection of $\sf{B}$ and $\sf{F}_n$ then they cannot be comparable. Hence we see that the elements of 
$\sf{B}  \cap \sf{F}_n$ are the characteristic vectors of a Sperner family in $[n]$. Let $m_n$ be the largest cardinality of a Sperner family of subsets of 
$[n]$ whose characteristic vectors are in $\sf{F}_n$. If we drop the last coordinate of the characteristic vectors, these remain distinct because if two vectors are 
incomparable, then they differ in at least two coordinates. Further, the shortened strings of length $n-1$ clearly belong to $\sf{F}_{n-1}$. This yields
$$m_n\leq f_{n-1}\leq \beta f_n$$
for some constant $\beta <1$ and every natural $n$, where the last inequality follows from the monotonicity and the well--known asymptotics of the standard Fibonacci sequence, according to which 
$\frac {f_{n-1}}{f_n}$ converges to $\frac{2}{1+\sqrt{5}}<1$.
Observing that $f_n\geq 2^{0.694n}$ for suffficiently large $n$, we conclude that 
$$\vert \overline{\sf{B} } \vert \geq f_n-m_n \geq (1-\beta) f_n\geq (1-\beta)2^{0.694n}>2^{0.69n}$$
for $n$ large enough. Hence 
$$\vert {\sf B} \vert \leq  2^n-2^{0.69n}$$
for sufficiently large $n$, as claimed.

To prove the lower bound we shall exhibit a set of pairwise skewincident sequences. The weight $w(\x)$ of a binary string $\x \in \{0,1\}^n$ is its number of 1's. In case of 
$\x=x_1x_2 \dots x_n$ we have 
$$w(\x):=\sum_{i=1}^n x_i$$
The support set of a string $\x \in \{0,1\}^n$ is the set ${\sf S}(\x)\subseteq [n]$ of positions $i$ in which $x_i=1$. In other words,
$w(\x)=\vert {\sf S}(\x) \vert$.
The influence ${\bf i}(\x)$ of string $\x$ is a binary string of the same length that has a 1 in position $j\in [n]$ if and only if either $x_{j-1}=1$ and/or 
$x_{j+1}=1$. We write
$$\gamma(\x):=w(\x)+w({\bf i}(\x))$$
and define the set ${\sf C}_n \subseteq  \{0,1\}^n$ as
$${\sf C}_n:=\left \{\x\, \vert \, \gamma(\x)>n \right \}$$

We claim that any two distinct elements of ${\sf C}_n$ are skew coincident. In fact, consider $\x \in {\sf C}_n$ and $\y \in {\sf C}_n$. Then we have 
$$w(\x) +w({\bf i}(\x))>n \;\hbox{and}\; w(\y) +w({\bf i}(\y)) >n$$
whence
\begin{equation}\label{eq: fir}
w(\x) +w({\bf i}(\x))+w(\y) +w({\bf i}(\y))>2n
\end{equation}
If $\x$ and $\y$ were not skew coincident, the sets ${\sf S}(\x)$ and ${\sf S}({\bf i}(\y))$ would be disjoint, implying that
$$w(\x)+w({\bf i}(\y)) \leq n$$
and likewise,
$$w(\y)+w({\bf i}(\x)) \leq n$$
yielding
$$w(\x)+w({\bf i}(\y))+w(\y)+w({\bf i}(\x))\leq 2n$$
in contradiction with (\ref{eq: fir}).

To lower bound the cardinality of ${\sf C}_n$ we shall use a well--known concentration inequality of McDiarmid \cite{M}. 
Let the random variable $X^n=X_1X_2\dots X_n$ be uniformly distributed on $\{0,1\}^n$. Then the variables $X_i$, $i\in [n]$ are 
totally independent and uniformly distributed over $\{0,1\}$. To prove our lower bound, it suffices to show that
\begin{equation}\label{mai}
{\rm Pr}\{\gamma(X^n)\leq n \} \leq 2^{-0.04n}
\end{equation}
Let $\alpha_i(\x)$ denote the $i$'th coordinate of the vector ${\bf i}(\x)$. We write
$$\gamma_i(\x):=x_i+\alpha_i(\x)$$
Hence 
$$\gamma(X^n)=\sum_{i=1}^n \gamma_i(X^n)$$
The function $\gamma(\x)$ defined on $ \{0,1\}^n$ satisfies the Lipschitz condition that given any two arguments $\x$ and $\y$ 
differing only in the $i$'th coordinate we have, for every $i \in [n]$
$$\vert \gamma(\x)-\gamma(\y) \vert \leq 2$$
This is simply because the possible values of $\gamma_i$ are only 0,1 or 2.
Let us now calculate the expected value of the random variable $\gamma(X^n)$.
By the linearity of the expected value and the definition of $\gamma_i$ we have
\begin{equation}\label{eq: ex}
\E \gamma(X^n)=\sum_{i=1}^n \E \gamma_i(X^n)=\sum_{i=1}^n[\E X_i+\E \alpha_i(X^n)]
\end{equation}
Since, for every $i$, both $X_i$ and $\alpha_i(X^n)$ take only the values 0 and 1, we have that
$$\E X_i+\E \alpha_i(X^n)={\rm Pr}\{X_i=1\}+{\rm Pr}\{\alpha_i(X^n)=1\}$$
Since $X_i$ is uniformly distributed, for every $i \in [n]$ 
$${\rm Pr}\{X_i=1\}=\frac{1}{2}$$
Also, since the $X_i$ are totally independent, and because $\alpha_i(X^n)=0$ if and only if $X_{i-1}=X_{i+1}=0$, 
for $1<i<n$ we see that 
$${\rm Pr}\{\alpha_i(X^n)=1\}=\frac{3}{4}$$
while ${\rm Pr}\{\alpha_i(X^n)=1\}=\frac{1}{2}$ else.
Thus we obtain
$$\E \gamma(X^n)= \frac{5n}{4}-\frac{1}{2}$$
and
$${\rm Pr}\{X^n \in  \overline{{\sf C}_n}\}\leq {\rm Pr}\{\gamma(X^n)\leq n\}\leq  
{\rm Pr}\left\{\vert \gamma(X^n)-\E (\gamma(X^n)) \vert>\frac{n}{4}-\frac{1}{2}  \right \}    $$
Upper bounding the right--most probability by (13) in Theorem 3.1 of McDiarmid \cite{M} we see that for large enough $n$
$${\rm Pr}\{X^n \in  \overline{{\sf C}_n}\}\leq \exp \left (-\frac{2n^2}{(\ln 2)\cdot 64 n}\right )\leq \exp( -0.04n)$$
\hfill$\Box$

\noindent{\bf Remark}

It is easy to see that the set of strings used to establish the lower bound does not have maximum cardinality. In fact, it is not even maximal.

\section{Generalizations}
The question about skewincidence can be generalized to a problem about subgraphs of an arbitrary finite graph. We will say that two subsets of the vertex set of a graph are 
neighbors if they contain two respective vertices that are adjacent in the graph. Note that a subset may or may not be its own neighbor. Let us denote by $M(G)$ the maximum number of distinct subsets of the vertex set of the graph such that any two of them are neighbors. For many graphs we will be able to completely determine this number. In particular, this is the case for complete bipartite graphs. Complete multi--partite graphs are equally easy to treat so that we omit the details. In case of other graphs things can be much more complicated. In particular,  it is easy to see that $M(n)=M(P_n)$ where $P_n$ is the path of $n$ vertices. 
In what follows, a stable set in a graph is a set of pairwise non--adjacent vertices. 
\begin{prop}\label{prop: bic}
Let $K_{m,n}$ be the bipartite complete graph whose maximal stable(edge--free) sets have $m$ and $n$ vertices, respectively. 
Then
$$M(K_{m,n})=(2^m-1)(2^n-1)+2$$
More generally, if $K_{n_1,n_2,\dots n_r}$ is a complete multipartite graph with disjoint stable sets of cardinality $n_1, n_2, \dots n_r$, respectively. 
We have 
$$M(K_{n_1,n_2,\dots n_r})=2^{\sum_{i=1}^r n_i}-\sum_{i=1}^r 2^{n_i}+2r-1$$
\end{prop}

\proof

It is obvious in the bipartite case that a family of subsets with the desired property cannot contain more than one subset of any of the two maximal stable sets. In the $r$--partite case for $r>2$, exactly in the same way, a family as required cannot contain more than one subset of any of the maximal stable sets.
\hfill$\Box$

All the above can be considered as special cases of a single more general problem other special cases of which contain the original Shannon set-up of graph capacity.

Let $F$ be a graph with vertex set $\N$ and $G$ arbitrary, finite or infinite. 
Consider, for every $n\in \N$ the family of all the mappings $f:[n]\rightarrow V(G)$ and denote it by $\M(F, G, n)$. We will say that 
two of these, $a \in \M(F, G, n)$, $b\in \M(F, G, n)$ form an attractive couple if there exist two, not necessarily distinct numbers $i\in [n]$ and 
$j\in [n]$ 
such that $i$ and $j$ are adjacent in $F$ while $a(i)$ and $b(j)$ are adjacent in $G$. We are interested in determining the largest cardinality of a subset of pairwise attractive elements of $\M(F, G, n)$.

If $F$ is the all--loops graph and $G$ an arbitrary simple graph, then $|\M(F, G, n)|$ is exponential in $n$ and the (always existing) limit of 
$\sqrt[n]{|\M(F, G, n)|}$  is the Shannon capacity of the graph $G$. If $F$ is the semi--infinite path and $G$ is a graph with two vertices and a loop as its only edge, we get back the problem of skew--incidence. Its immediate generalizations are obtained if $F$ is arbitrary while $G$ remains the same one--edge graph as for the skew--incidence problem.

If $G$ also has $\N$ as its vertex set then we will sometimes restrict attention to the subset $\B(M, G, n)\subseteq \M(F, G, n)$ of bijective mappings 
from $[n]$ onto itself. This leads, in case of the all--loop graph in the role of $F$ to the concept of permutation capacity. 

\section{A Sperner--type problem}

As a byproduct from the proof of the Theorem, we get the following extremely simple sounding problem in classical extremal set theory. Let $\sf{F}_n$ be the 
set of all the binary sequences of length $n$ without 1's in consecutive positions. (Their numbers $f_n$ are the classical example for the standard Fibonacci sequence.)
We consider these binary sequences as the characteristic vectors of subsets of the set $[n]$ in the usual manner and ask for the maximum cardinality of a Sperner family they
contain. 

In our proof a very weak upper bound on this cardinality was sufficient. The present problem is interesting inasmuch no classical proof for Sperner's theorem \cite{Sp}
seems to be suitable to solve it.

\end{document}